\documentclass[12pt]{article}

\usepackage{xypic,amsthm,amssymb,amsmath}

\newtheorem{theorem}{Theorem}
\newtheorem{lemma}[theorem]{Lemma}
\newtheorem{corollary}[theorem]{Corollary}
\newtheorem{proposition}[theorem]{Proposition}

\theoremstyle{remark}
\newtheorem{remark}[theorem]{Remark}

\allowdisplaybreaks

\newcommand{\divides}{\,|\,}
\newcommand{\tensor}{\otimes}
\newcommand{\lra}{\longrightarrow}

\newcommand{\setofall}[2]{\mbox{$\{\,#1\,|\,#2\,\}$}}

\newcommand{\Sym}{\mathcal{S}}
\newcommand{\CC}{\mathcal{C}}
\newcommand{\DD}{\mathcal{D}}

\newcommand{\Z}{\mathbb{Z}}
\newcommand{\Q}{\mathbb{Q}}

\newcommand{\abacus}[2]{\langle\,#1,#2\,\rangle}
\DeclareMathOperator{\Rad}{rad}
\DeclareMathOperator{\ch}{char}
\DeclareMathOperator{\Hom}{Hom}
\DeclareMathOperator{\GL}{GL}
\DeclareMathOperator{\sgn}{sgn}

\begin{document}

\title{Modular Lie Powers\\ and the Solomon descent algebra}

\author{%
  \begin{minipage}[t]{58mm}
    \begin{center}
      Karin Erdmann\\
    \texttt{erdmann@maths.ox.ac.uk}
    \end{center}
  \end{minipage}
  \hspace*{\fill}
  and
  \hspace*{\fill}
  \begin{minipage}[t]{58mm}
    \begin{center}
      Manfred Schocker\thanks{%
        supported by Deutsche Forschungsgemeinschaft (DFG-Scho 799)}\\
      \texttt{schocker@maths.ox.ac.uk}
    \end{center}
  \end{minipage}\\[10mm]
  Mathematical Institute\\[1mm]
  24--29 St. Giles'\\[1mm]
  Oxford OX1 3LB}

\date{}

\maketitle

\begin{center}
  MSC 2000: 17B60 (primary), 17B01, 20C20, 20C30, 05E99 (secondary)
\end{center}

\newpage

\begin{abstract} \noindent
  Let $V$ be an $r$-dimensional vector space over an infinite field
  $F$ of prime characteristic $p$, and let $L_n(V)$ denote the $n$-th
  homogeneous component of the free Lie algebra on $V$. We study the
  structure of $L_n(V)$ as a module for the general linear group
  $GL_r(F)$ when $n=pk$ and $k$ is not divisible by $p$ and where $n
  \geq r$.  Our main result is an explicit 1-1 correspondence,
  multiplicity-preserving, between the indecomposable direct summands of
  $L_k(V)$ and the indecomposable direct summands of $L_n(V)$ which are not
  isomorphic to direct summands of $V^{\otimes n}$.  The direct
  summands of $L_k(V)$ have been parametrised earlier, by Donkin and
  Erdmann.  Bryant and St\"{o}hr have considered the case $n=p$ but
  from a different perspective.\\  
  Our approach uses idempotents of the Solomon descent algebras, and
  in addition a correspondence theorem for permutation modules of
  symmetric groups.
\end{abstract}

\newpage

\section{Introduction}

Let $F$ be an infinite field of prime characteristic $p$ and let $V$
be an $r$-dimensional vector space over $F$.  Let $L(V)$ be the free
Lie algebra on $V$ and denote its homogeneous component of degree $n$
by $L_n(V)$, for each positive integer $n$.  The group of graded
automorphisms of $L(V)$ can be identified with the general linear
group $\GL_r(F)$ in such a way that $L_1(V)$ becomes the natural
$\GL_r(F)$-module. In this way, $L_n(V)$ becomes a submodule of the
$n$-fold tensor product $V^{\tensor n}$, called the $n$th Lie power of
$V$.  One would like to know the structure of $L_n(V)$ as a module for
$\GL_r(F)$. As is well-known, if $p$ does not
divide $n$ then $L_n(V)$ is a direct summand of $V^{\tensor n}$. The
direct sum decomposition of $L_n(V)$ when $p$ does not divide $n$ was
dealt with in \cite{donkin-erdmann98}, generalising naturally the
classical theorems for characteristic zero.  Degrees divisible by $p$,
however, have largely been a mystery.

Here we study the module structure of $L_n(V)$ when $n$ is divisible
by $p$ but not divisible by $p^2$.  Our main result is the following:

\begin{theorem} \label{GL-haupt}
  Let $n=pk$ such that $k$ is not divisible by $p$ and assume $r\ge
  n$.  Then there is a 1-1 correspondence, multiplicity-preserving,
  between the indecomposable direct summands of $L_k(V)$ and the
  indecomposable direct summands of $L_n(V)$ which are not isomorphic
  to direct summands of $V^{\tensor n}$.
\end{theorem}

The case $k=1$ was considered in~\cite{bryant-stoehr-pre}.
 
The Schur functor relates representations of $\GL_r(F)$ with
representations of the symmetric group $\Sym_n$ for $r\ge n$
(see~\cite[Chapter~6]{green80}).  The image under the Schur functor of
$L_n(V)$ is the Lie module $L_n$ of $\Sym_n$ over $F$ which can be
described as $L_n=\omega_nF\Sym_n$ where
$$
\omega_n
=
(1-\zeta_n)(1-\zeta_{n-1})\cdots (1-\zeta_2)\in F\Sym_n
$$
is the Dynkin operator (see, for instance,~\cite{blelau93a}). Here
$\zeta_k$ denotes the descending $k$-cycle $(k\,\ldots\,1)\in \Sym_n$.

Our main tools come from the Solomon descent algebra, a subalgebra of
the group algebra of the symmetric group which contains $\omega_n$
(see Section~\ref{2}).  An analogous approach was used
in~\cite{schocker01} to study Lie powers over fields of characteristic
zero.

Let $S^p(L_k)$ be the image under the Schur functor of the $p$th
symmetric power $S^p(L_k(V))$ of $L_k(V)$.  In Section~\ref{4}, we
shall show that there is a short exact sequence of $\Sym_n$-modules
$$
0 \to L_n \to e_nF{\cal S}_n \to S^p(L_k) \to 0
$$
where $e_n$ is an idempotent in the Solomon descent algebra.

The middle term is projective, hence the Heller operator $\Omega$
gives us a {1-1} correspondence between the non-projective
indecomposable direct summands of $S^p(L_k)$ and the non-projective
indecomposable direct summands of $L_n$.  The module $S^p(L_k)$ turns
out to be a direct summand of a permutation module and can be analysed
with modular representation theory tools (see Section~\ref{5}). As a
result, there is a {1-1} correspondence, multiplicity-preserving,
between the non-projective indecomposable direct summands of
$S^p(L_k)$ and the indecomposable direct summands of $L_k$.  (In fact,
this result holds more generally, see Theorem~\ref{param}.)

This implies:

\begin{theorem} \label{sym-haupt}
  Let $k$ be a positive integer not divisible by $p$, then there is a
  {1-1} correspondence, multiplicity-preserving, between the
  non-projective indecomposable direct summands of $L_{pk}$ and the
  indecomposable direct summands of $L_k$.
\end{theorem}

For instance, if $p$ is odd then $L_{2p}$ has a unique non-projective
indecomposable direct summand, since $L_2$ has dimension one. A detailed
description of this non-projective summand of $L_{2p}$ is given in
Theorem~\ref{2p}.

The short exact sequence of $\Sym_n$-modules has an analogue on the
level of $\GL_r(F)$-modules, namely
$$
0 \to L_n(V) \to e_nV^{\tensor n} \to S^p(L_k(V)) \to 0.
$$
This is shown in Section~\ref{6}.  The modules occuring in this
sequence are $n$-homogeneous polynomial representations of $\GL_r(F)$,
that is, they are modules for the Schur algebra $S(r,n)$,
see~\cite[Chapter~2]{green80}.  For $r\ge n$, the indecomposable
direct summands of $V^{\tensor n}$ are precisely the indecomposable
modules which are projective and injective as modules for the Schur
algebra \cite[p.~94]{donkin98}.  The middle term, a direct summand of
$V^{\tensor n}$, is thus projective and injective as a module for
$S(r,n)$.  It follows that $\Omega$ gives a {1-1} correspondence
between the indecomposable direct summands of $S^p(L_k(V))$ which are
not projective and injective and the indecomposable direct summands of
$L_n(V)$ which are not projective and injective.  A {1-1}
correspondence between the indecomposable direct summands of
$S^p(L_k(V))$ which are not projective and injective and the
indecomposable direct summands of $L_k(V)$ is then readily derived and
allows us to deduce Theorem~\ref{GL-haupt} (see
Proposition~\ref{GL-param} and the subsequent remark).

As a further consequence, every indecomposable direct summand of $L_{pk}$ is
liftable, hence the formal character of any indecomposable direct summand of
$L_{pk}(V)$ is a sum of Schur functions (see
Corollary~\ref{liftable}).

Using the short exact sequence of $\GL(V)$-modules mentioned above, we
also prove a special case of a conjecture of Bryant~\cite{bryant-2} on
Lie resolvents in the final Section~\ref{7}.

In what follows we shall simply say ``summand'' for ``direct summand'' of
a module.

\section{The descent algebra} \label{2}

Let $n$ be a positive integer.  We summarise properties of the descent
algebra of the symmetric group $\Sym_n$. For general reference,
see~\cite{blelau93a,garsia-reutenauer89,reutenauer93,schocker-fields}.
Note that products $\pi\sigma$ of permutations $\pi,\sigma\in\Sym_n$
are to be read from left to right: first $\pi$, then $\sigma$.

Let $\mu$ be a composition of $n$, that is, a finite sequence
$(\mu_1,\ldots,\mu_k)$ of positive integers with sum $n$.  We then
write $\Sym_\mu$ for the usual embedding of the direct product
$\Sym_{\mu_1}\times\cdots\times \Sym_{\mu_k}$ in $\Sym_n$.

The length of a permutation $\pi$ in $\Sym_n$ is the number of
inversions of $\pi$. Each right coset of $\Sym_\mu$ in $\Sym_n$
contains a unique permutation of minimal length. Define $X^\mu$ to be
the sum in the integral group ring $\Z \Sym_n$ of all these minimal
coset representatives of $\Sym_\mu$ in $\Sym_n$.  For example,
$X^{(n)}$ is the identity of $\Sym_n$, while $X^{(1,1\ldots,1)}$ is
the sum over all permutations in $\Sym_n$.  Due to Solomon
\cite[Theorem~1]{solomon76}, the $\Z$-linear span $\DD_n$ of the
elements $X^\mu$ ($\mu$ any composition of $n$), is a subring of $\Z
\Sym_n$ of rank $2^{n-1}$, called the \emph{descent algebra}
of~$\Sym_n$. In fact, the elements $X^\mu$ form a $\Z$-basis of
$\DD_n$ and there exist nonnegative integers $c_{\lambda\mu\nu}$ such
that
\begin{equation}  \label{solorule}
  X^\lambda\,X^\mu=\sum_\nu c_{\lambda\mu\nu}\,X^\nu,
\end{equation}
for all compositions $\lambda$, $\mu$ of $n$. It is well-known that
$\omega_n\in\DD_n$ and a result of Dynkin-Specht-Wever states that
$\omega_n^2=n\omega_n$ (see, for instance, \cite{blelau93a,garsia89}).

The \emph{Young character} $\varphi^\mu$ of $\Sym_n$ is induced from
the trivial character of $\Sym_\mu$; that is, $\varphi^\mu(\pi)$ is
the number of right cosets of $\Sym_\mu$ in $\Sym_n$ which are fixed
by right multiplication with $\pi$, for any $\pi\in \Sym_n$. The
$\Z$-linear span $\CC_n$ of the Young characters $\varphi^\mu$ ($\mu$
any composition of $n$) is a subring of the ring of $\Z$-valued class
functions of $\Sym_n$.  In fact, the elements $\varphi^\mu$ form a
$\Z$-basis of $\CC_n$ and products have the form
\begin{equation}  \label{youngrule}
  \varphi^\lambda\,\varphi^\mu=\sum_\nu c_{\lambda\mu\nu}\,\varphi^\nu
\end{equation}
for all compositions $\lambda$, $\mu$ of $n$, with the same
coefficients as in~\eqref{solorule}. As a consequence, the $\Z$-linear
map $c_n:\DD_n\to \CC_n$, defined by
\begin{equation}  \label{solhom}
  X^\mu\longmapsto \varphi^\mu
\end{equation}
for all compositions $\mu$ of $n$, is an epimorphism of rings.  This
is the second part of \cite[Theorem~1]{solomon76}.

\begin{theorem}
  Let $F$ be a field, then the $F$-linear span $\DD_{n,F}$ of the
  elements $X^\mu$ is a subalgebra of the group algebra $F\Sym_n$, while
  the $F$-linear span $\CC_{n,F}$ of the $F$-valued Young characters
  $$
  \varphi^{\mu,F}:\Sym_n\to F,\,\pi\mapsto\varphi^\mu(\pi)\cdot 1_F
  $$
  is a subalgebra of the algebra of $F$-valued class functions of
  $\Sym_n$.  
  The $F$-linear map 
  $$
  c_{n,F}:\DD_{n,F}\to \CC_{n,F}
  $$
  sending $X^\mu$ to $\varphi^{\mu,F}$ for all compositions $\mu$
  of $n$, is an epimorphism of algebras.
\end{theorem}

Indeed, by definition, there are the product formulae~\eqref{solorule}
and~\eqref{youngrule} in $\DD_{n,F}$ and $\CC_{n,F}$, respectively,
where the coefficients $c_{\lambda\mu\nu}$ should be read as
$c_{\lambda\mu\nu}\cdot 1_F$.

\begin{theorem} \label{atwill}
  $\Rad\,\DD_{n,F}=\ker\,c_{n,F}$.
\end{theorem}

This is \cite[Theorem~3]{solomon76} if $F$ has characteristic zero and
\cite[Theorem~2]{aw97} if $F$ has prime characteristic.  As a
consequence, $\DD_{n,F}/\Rad\,\DD_{n,F}\cong\CC_{n,F}$.

We analyse the algebra $\CC_{n,F}$. The conjugacy classes $C_\lambda$
of $\Sym_n$ are indexed by partitions $\lambda$ of $n$, in a natural
way: $C_\lambda$ consists of all permutations in $\Sym_n$ of cycle
type $\lambda$.

Let $\lambda$, $\mu$ be partitions of $n$ and $p$ be a prime, then
$\lambda$ and $\mu$ are \emph{$p$-equivalent} if the $p$-regular parts
of $\pi$, $\sigma$ are conjugate in $\Sym_n$, for each $\pi\in
C_\lambda$, $\sigma\in C_\mu$.  Note that the cycle type $\nu$ of the
$p$-regular part of $\pi\in C_\lambda$ is obtained from $\lambda$ by
replacing each entry $\lambda_i=kp^m$ of $\lambda$ by the sequence
$(k,\ldots,k)$ of length $p^m$, where $k\ge 1$ and $m\ge 0$ are so
chosen that $k$ is not divisible by $p$. For example, if $p=2$ and
$\lambda=(6,3,2)$, then $\nu=(3,3,3,1,1)$.

The partition $\mu$ is \emph{$p$-regular} if no part of $\mu$ occurs
more than $p-1$ times in $\mu$.  It is convenient to extend these
definitions to the case where $p=0$, by saying that any partition is
$0$-regular, and $0$-equivalent to itself only.

Then the $p$-regular partitions form a transversal for the
$p$-equivalence classes of partitions, so that each partition
$\lambda$ is $p$-equivalent to a unique $p$-regular partition $\mu$.

For the remainder of this section, $F$ is a field of characteristic
$p$ (which might be zero or not).  Define $C_{\mu,F}$ to be the union
of all conjugacy classes $C_\lambda$ in $\Sym_n$ such that $\lambda$ is
$p$-equivalent to $\mu$, for each $p$-regular partition $\mu$ of $n$.
Let $\ch_{\mu,F}$ denote the characteristic function $\Sym_n\to F$ of
$C_{\mu,F}$ (mapping $\pi\in \Sym_n$ to $1_F$ or zero according as
$\pi\in C_{\mu,F}$ or not).

\begin{proposition} \label{apw}
  $\CC_{n,F}$ is split semisimple. The elements $\ch_{\mu,F}$, indexed
  by $p$-regular partitions $\mu$ of $n$, form a full set of primitive
  idempotents in $\CC_{n,F}$.
\end{proposition}
(see, for instance, \cite[Lemma~2 and its proof]{aw97}) Combining
Proposition~\ref{apw} with Theorem~\ref{atwill}, we obtain the
following result.

\begin{corollary} \label{pri-id-d}
  $\DD_{n,F}$ is a basic algebra with irreducible modules indexed by
  $p$-regular partitions of $n$. In fact, there exists a complete set
  of mutually orthogonal primitive idempotents
  $$
  \setofall{e_{\mu,F}}{\mbox{$\mu$ $p$-regular}}
  $$
  in $\DD_{n,F}$ such that $c_{n,F}(e_{\mu,F})=\ch_{\mu,F}$ for all
  $\mu$.
\end{corollary}

(for the second part, see \cite[Theorem~44.3]{dornhoff72B}.)

The idempotent $e_n=e_{(n),F}$ will be of crucial importance for our
study of the modular Lie representations.  This section concludes with
two observations on $e_n$.  Choose coefficients $a_\mu\in F$ such that
$e_n=\sum_\mu a_\mu X^\mu$, where the sum is over compositions $\mu$
of $n$. The image of $e_n$ is $\ch_{(n),F}$, thus maps long cycles
$\pi\in C_{(n)}$ to~$1_F$. But $\varphi^{\mu,F}(\pi)=0$ for all such
$\pi$ whenever $\mu\neq(n)$, while $\varphi^{(n),F}(\pi)=1_F$.  This
implies $a_{(n)}=1_F$, that is
\begin{equation}  \label{1-1}
  e_n=X^{(n)}+\sum_{\mu\neq(n)} a_\mu X^\mu.
\end{equation}
The second important property of the idempotent $e_n$ is that
\begin{equation}  \label{1-2}
  \dim\,e_n F\Sym_n
  =
  |C_{(n),F}|.
\end{equation}
This is a special case of the following result.

\begin{proposition}
  $\dim\,e_{\mu,F} F\Sym_n=|C_{\mu,F}|$, for each $p$-regular
  partition $\mu$ of $n$.
\end{proposition}

\begin{proof}  
  Let $M_{\mu,F}$ denote the (one-dimensional) irreducible
  $\DD_{n,F}$-module corresponding to $e_{\mu,F}$, for each
  $p$-regular partition $\mu$ of $n$.  The action of
  $\alpha\in\DD_{n,F}$ on $M_{\mu,F}$ is then scalar multiplication
  with $c_{n,F}(\alpha)(\pi)$, where $\pi\in C_\mu$.  In particular,
  the family $\{M_{\mu,F}\}$ is defined over $\Z$, since the Young
  characters take values in $\Z$.
  
  The decomposition matrix of $\DD_{n,\Q}$ modulo $p$ is very simple;
  if $\lambda$ and $\mu$ are partitions of $n$ and the characteristic
  $p$ of $F$ is positive, then $M_{\lambda,F}\cong
  M_{\lambda,\Z}\tensor F$ and $M_{\mu,F}\cong M_{\mu,\Z}\tensor F$
  are isomorphic as $\DD_{n,F}$-modules if and only if $\lambda$ and
  $\mu$ are $p$-equivalent. This is due to
  Atkinson-Pfeiffer-Willigenburg \cite[Theorem~4 and
  Section~4.1]{apw02} and follows from Theorem~\ref{atwill} and
  Proposition~\ref{apw}.
    
  If $V$ is an arbitrary $\DD_{n,F}$-left module, then the
  multiplicity $[V:M_{\mu,F}]$ of $M_{\mu,F}$ in a composition series
  of $V$ is equal to the dimension of $e_{\mu,F}V$, since $M_{\mu,F}$
  has dimension one. It follows that
  $$
  \dim\,e_{\mu,F} F\Sym_n
  =
  [F\Sym_n:M_{\mu,F}]
  =
  \sum_\lambda [\Q \Sym_n:M_{\lambda,\Q}]
  =
  \sum_\lambda \dim\,e_{\lambda,\Q}\Q \Sym_n\,,
  $$
  where both sums are taken over all partitions $\lambda$ of $n$
  which are $p$-equivalent to~$\mu$.  However, $\dim\,e_{\lambda,\Q}\Q
  \Sym_n=|C_\lambda|$ for all partitions $\lambda$ of $n$ (see, for
  instance, \cite[Lemma~3.2]{schocker-fields}). This completes the
  proof.
\end{proof}

\section{The Lie module in prime degree}  \label{3}

Let $F$ be a field of prime characteristic $p$. To illustrate our
approach, we start with studying the Lie module $L_p$ of the symmetric
group $\Sym_p$ (although this will be generalised later).  This module
has already been analysed in~\cite{bryant-stoehr-pre,bks99}.

For each positive integer $n$, we consider the sum $s_n$ in $F\Sym_n$
of all permutations in $\Sym_n$ as a linear generator of the trivial
$\Sym_n$-module $F$.

\begin{theorem} \label{p}
  There is a short exact
  sequence of $\Sym_p$-right modules
  $$
  0
  \lra
  L_p
  \stackrel{\alpha}{\lra} 
  e_p F\Sym_p
  \stackrel{\beta}{\lra} 
  F
  \lra 
  0,
  $$
  where $\alpha$ is left multiplication with $e_p$ and $\beta$ is
  left multiplication with $s_p$.
\end{theorem}

\begin{proof}
  First let $n$ be an arbitrary positive integer, then a
  multiplication rule in $\DD_n$ (over $\Z$) is
  \begin{equation}  \label{xiom}
    \mbox{$X^\mu\omega_n=0$ whenever $\mu\neq(n)$}
  \end{equation}
  (see \cite[\S 2]{garsia89}), which implies
  \begin{equation}  \label{enon}
    e_n\omega_n
    =
    \omega_n+\sum_{\mu\neq(n)} a_\mu X^\mu\omega_n
    =
    \omega_n\,,
  \end{equation}
  by \eqref{1-1}.  Thus, in particular, $\alpha$ is an inclusion.
  Furthermore, the image $\omega_pF\Sym_p$ of $\alpha$ is contained in
  the kernel of $\beta$, since
  $s_p\omega_p=X^{(1,\ldots,1)}\omega_p=0$, by \eqref{xiom}.
  
  If $\mu=(\mu_1,\ldots,\mu_l)$ is a composition of $p$, then the
  number of summands in $X^\mu$ is equal to the number of right cosets
  of $\Sym_\mu$ in $\Sym_p$ which is
  $\binom{p}{\mu_1\;\cdots\;\mu_l}$. Hence $s_pX^\mu=0$ in $F\Sym_p$
  whenever $\mu\neq(p)$.  This implies $ s_p e_p = s_p $, by
  \eqref{1-1} again. Thus $\beta$ is onto.
  
  Finally, $\dim\,L_p=(p-1)!$ is well-known, while
  $\dim\,e_pF\Sym_p=(p-1)!+1$ follows from \eqref{1-2}.  Comparing
  dimensions, completes the proof.
\end{proof}

\begin{corollary}
  $L_p$ has a unique non-projective indecomposable summand, which is
  isomorphic to the Specht module associated to the partition
  $(p-1,1)$.
\end{corollary}

\begin{proof}
  The module $e_pF\Sym_p$ is projective, since $e_p$ is an idempotent.
  The trivial $\Sym_p$-module is non-projective indecomposable, hence
  by Theorem~\ref{p}, the only non-projective indecomposable summand
  of $L_p$ is isomorphic to the Heller translate $\Omega(F)$ of $F$
  (see~\cite[\S 1.5]{benson98}). This is known to be the Specht module
  mentioned.
\end{proof}

The character of $L_p$ over $\Q$ is known, and also knowing that $L_p$
is the direct sum of a certain Specht module and a projective module,
determines uniquely the projective summand. Details have been worked
out in \cite[Theorem~6.2]{bryant-stoehr-pre}; see also~\cite{bks99}.

\section{Lie modules in degree not divisible by $p^2$}  \label{4}

With a little more input from the descent algebra we will now extend
Theorem~\ref{p} to the Lie module $L_n$ of $\Sym_n$ for arbitrary $n$
divisible by $p$, but not divisible by $p^2$.

Throughout, $F$ is a field of prime characteristic $p$ and $n=kp$ with
$k$ not divisible by $p$. The \emph{$p$th symmetrisation} $S^p(L_k)$
of the Lie module for $\Sym_k$ is a module for $\Sym_n$ and is defined
as a module induced from the wreath product $\Sym_k\wr \Sym_p$, as
follows.

For $\alpha\in \Sym_a$, $\beta\in \Sym_b$, define $\alpha\#\beta\in
\Sym_{(a,b)}\subseteq \Sym_{a+b}$ in the natural way: $\alpha\#\beta$ maps
$i$ to $i\alpha$ if $i\le a$, and to $(i-a)\beta+a$ otherwise. If,
additionally, $\gamma\in \Sym_c$, then
$(\alpha\#\beta)\#\gamma=\alpha\#(\beta\#\gamma)$, as is readily seen.
Using linearity, we define
$$
\omega^{(k,\ldots,k)}
:=
\omega_k\#\cdots\#\omega_k
\in
F(\Sym_k\#\cdots\# \Sym_k)
=
F\Sym_{(k,\ldots,k)}
\quad
\mbox{($p$ factors)}.
$$

Now assume $\pi\in \Sym_p$. Let $\pi^{[k]}$ be the element of
$\Sym_{kp}$ which permutes the $p$ successive blocks of size $k$ in
$\{1,\ldots,kp\}$ according to $\pi$. More explicitly, set
$(ik-j)\pi^{[k]}=(i\pi)k-j$ for all $i\in\{1,\ldots,p\}$ and
$j\in\{0,\ldots,k-1\}$.

The map $\pi\mapsto \pi^{[k]}$ extends linearly to an embedding of
$F\Sym_p$ in $F\Sym_{kp}$, since
$(\pi\sigma)^{[k]}=\pi^{[k]}\sigma^{[k]}$ for all $\pi,\sigma\in
\Sym_k$.  Furthermore, $\Sym_p^{[k]}\Sym_{(k,\ldots,k)}$ is isomorphic
to the wreath product $\Sym_k\wr \Sym_p$, since
\begin{equation} \label{wr-def1}
  (\alpha_1\#\cdots\#\alpha_p)(\beta_1\#\cdots\#\beta_p)
  =
  (\alpha_1\beta_1)\#\cdots\#(\alpha_p\beta_p)
\end{equation}
and
\begin{equation} \label{wr-def2}
  \pi^{[k]}(\alpha_1\#\cdots\#\alpha_p)
  =
  (\alpha_{1\pi}\#\cdots\#\alpha_{p\pi})\pi^{[k]}
\end{equation}
for all $\pi\in\Sym_p$ and
$\alpha_1,\ldots,\alpha_p,\beta_1,\ldots,\beta_p\in\Sym_k$.

We define now the $p$th symmetrisation of $L_k$ by
$$
S^p(L_k)
=
s_p^{[k]}\omega^{(k,\ldots,k)}F\Sym_n\,.
$$
We shall see after Corollary~\ref{GL-kp} that this is isomorphic
to the image of the $p$th symmetric power $S^p(L_k(V))$ of $L_k(V)$
under the Schur functor.

The main result of this section is:

\begin{theorem} \label{kp}
  There is a short exact sequence of $\Sym_n$-right modules
  $$
  0
  \lra
  L_n
  \stackrel{\alpha}{\lra} 
  e_n F\Sym_n
  \stackrel{\beta}{\lra} 
  S^p(L_k)
  \lra 
  0,
  $$
  where $\alpha$ is left multiplication with $e_n$ and $\beta$ is
  left multiplication with $X^{(k,\ldots,k)}$.
\end{theorem}

This reduces to Theorem~\ref{p} in case $k=1$.  Note that
$\dim\,S^p(L_k)=|C_{(k,\ldots,k)}|$, so that
$$
\dim\,e_nF\Sym_n=|C_{(n)}|+|C_{(k,\ldots,k)}|=\dim\,L_n+\dim\,\Sym_p(L_k),
$$
by \eqref{1-2}.

From now on, we write $\kappa$ for the composition $(k,\ldots,k)$ of
$n=kp$. Before we prove Theorem~\ref{kp}, let us recall some more
multiplication rules for the descent algebra.

If $\lambda=(\lambda_1,\ldots,\lambda_l)$ and $\nu$ are compositions
of $n$, then write $\nu\le \lambda$ if there is a composition
$\nu^{(i)}$ of $\lambda_i$ for all $i\le l$ such that $\nu$ is equal
to the concatenation $(\nu^{(1)},\ldots,\nu^{(l)})$.  For example,
$(1,2,3,2,1,2)\le(3,3,5)$.

Concerning the coefficients $c_{\lambda\mu\nu}$ in \eqref{solorule},
there is the restriction
\begin{equation}  \label{restrict}
  \mbox{$c_{\lambda\mu\nu}=0$ unless $\nu\le \lambda$}
\end{equation}
\cite[Lemma~1(i)]{aw97}. Furthermore, for any composition
$\mu=(\mu_1,\ldots,\mu_l)$ of $n$ with each part divisible by $k$,
\begin{equation} \label{binomi}
  c_{\kappa\mu\kappa}=\binom{p}{\mu_1/k\cdots\mu_l/k}.
\end{equation}
This follows directly from the combinatorial description of
$c_{\lambda\mu\nu}$ given in \cite[Eq.~(1.1)]{aw97}, for example.
Finally,
\begin{equation}  \label{xiok}
  \mbox{$X^\mu\omega^\kappa=0$ unless $\kappa\le \mu$}
\end{equation}
and
\begin{equation} \label{gr21}
  X^\kappa\omega^\kappa=s_p^{[k]}\omega^\kappa
\end{equation}
(see \cite[Theorem~2.1]{garsia-reutenauer89}). It should be mentioned
that the equations~\eqref{xiok} and~\eqref{gr21} are much deeper than
the equations~\eqref{restrict} and~\eqref{binomi}.

\begin{proof}[Proof of Theorem~\ref{kp}]
  We have already seen that $e_n\omega_n=\omega_n$ (see~\eqref{enon}),
  so $\alpha$ is an inclusion.  Furthermore, $X^\kappa\omega_n=0$ as
  in the proof of Theorem~\ref{p}, by~\eqref{xiom}, so $L_n$ is
  contained in the kernel of $\beta$. It remains to show that
  $X^\kappa e_nF\Sym_n$ contains $S^p(L_k)$, for then a comparison of
  dimensions completes the proof.
  
  But, for each composition $\mu$ of $n$, \eqref{xiok} yields $X^\mu
  \omega^\kappa=0$ unless $\kappa\le\mu$ (that is, unless each part of
  $\mu$ is divisible by $k$).  In this case,
  \begin{eqnarray*}
    X^\kappa X^\mu \omega^\kappa
    & = &
    \sum_{\nu\le\kappa} c_{\kappa\mu\nu} X^\nu \omega^\kappa,
    \quad
    \mbox{ by~\eqref{solorule}, \eqref{restrict}}\\[1mm]
    & = &
    c_{\kappa\mu\kappa} X^\kappa \omega^\kappa,
    \quad
    \mbox{ by~\eqref{xiok}}\\[1mm]
    & = &
    \binom{p}{\mu_1/k\cdots\mu_l/k} X^\kappa \omega^\kappa,
    \quad
    \mbox{ by~\eqref{binomi}}.
  \end{eqnarray*}
  As a consequence, $X^\kappa e_n
  \omega^\kappa=X^\kappa\omega^\kappa=s_p^{[k]}\omega^\kappa$ in
  $F\Sym_n$, by~\eqref{1-1} and~\eqref{gr21}, which implies $X^\kappa e_n
  F\Sym_n\supseteq X^\kappa e_n\omega^\kappa F\Sym_n=s_p^{[k]}\omega^\kappa
  F\Sym_n=S^p(L_k)$, and we are done.
\end{proof}

As in the special case where $k=1$, Theorem~\ref{kp} gives the $1$-$1$
correspondence $U\mapsto \Omega(U)$ between the non-projective
indecomposable summands of $S^p(L_k)$ and those of $L_n$, since
$e_nF\Sym_n$ is a projective $\Sym_n$-module.

\section{Non-projective indecomposable summands\\ of symmetrised modules}
 \label{5}

Let $F$ be a field of prime characteristic $p$ and $k$ a positive
integer not divisible by $p$.  Bearing in mind Theorem~\ref{kp}, we
are aiming at a parametrisation of the non-projective indecomposable
summands of $S^p(L_k)$.

We shall consider, more generally, an arbitrary idempotent $e\in
F\Sym_k$ (instead of $\frac{1}{k}\,\omega_k$) and the corresponding
right ideal $U=eF\Sym_k$ of $F\Sym_k$. Let $n=kp$. Generalising the
definition of $S^p(L_k)$ in Section~\ref{4}, we let the $p$th
symmetrisation of $U$ be the $F\Sym_n$-module
$$
S^p(U)
=
s_p^{[k]}e^{\# p} F\Sym_n\,,
$$
where $e^{\# p}=e\#\cdots\# e$ ($p$ factors).
The aim is to prove

\begin{theorem} \label{param}
  There is a {1-1} correspondence, multiplicity-preserving, between
  the non-projective indecomposable summands of $S^p(U)$ and the
  indecomposable summands of $U$.
\end{theorem}

We write $H=\Sym_k\wr \Sym_p$ (which is taken as the subgroup
$\Sym_p^{[k]}\Sym_{(k,\ldots,k)}$ of $\Sym_n$).

A crucial step towards a proof of Theorem~\ref{param} is the following
observation.

\begin{proposition} \label{p-permu}
  We have 
  $S^p(U)\cong(F\tensor_{F\Sym_p} e^{\# p}FH)\tensor_{FH} F\Sym_n$.    
  In particular, $S^p(U)$ is a summand of a permutation module
  of $\Sym_n$.
\end{proposition}

\begin{proof}
  The element $e^{\# p}$ is an idempotent in $H$, by~\eqref{wr-def1},
  which commutes with every element of $\Sym_p\subseteq H$,
  by~\eqref{wr-def2}.  So $FH=e^{\# p}FH\oplus (1-e^{\# p})FH$ as a
  $(F\Sym_p,FH)$ bimodule. Hence $F\tensor_{F\Sym_p} e^{\# p} FH$ is a
  summand of $F\tensor_{F\Sym_p} FH$, which is a permutation module of
  $H$, and therefore $(F\tensor_{F\Sym_p}e^{\#
    p}FH)\tensor_{FH}F\Sym_n$ is a summand of a permutation module of
  $\Sym_n$. We want to identify this summand with $S^p(U)$.
  
  In general, for a subgroup $Y$ of a finite group $X$ and an
  idempotent $f$ of $FX$ which commutes with all elements of $Y$ we have
  $$
  gFY\tensor_{FY} fFX\cong gfFX
  $$
  for all $g\in FY$, since $fFX$ is a projective left $FY$-module.
  We apply this twice, the first time with $Y=\Sym_p$, $X=H$, $f=e^{\#
    p}$ and $g=s_p^{[k]}$ and the second time with $Y=H$, $X=\Sym_n$,
  $f=1$ and $g=s_p^{[k]}e^{\# p}$. This yields
  $$
  (F\tensor_{F\Sym_p} e^{\# p}FH)\tensor_{FH} F\Sym_n
  \cong
  s_p^{[k]}e^{\# p}FH\tensor_{FH} F\Sym_n
  \cong
  s_p^{[k]}e^{\# p}F\Sym_n
  =
  S^p(U)
  $$
  as asserted.
\end{proof}

For an arbitrary finite group $X$, an $X$-module $W$ is said to be a
\emph{$p$-permuta\-tion module} if, for every $p$-subgroup $P$ of $X$,
there is a $P$-invariant basis of $W$. Equivalently, $W$ is a 
summand of a permutation module of $X$.

The above proposition implies that $S^p(U)$ is a $p$-permutation
module for $\Sym_n$. In order to parametrise the non-projective
indecomposable summands of $S^p(U)$, it is natural to use their
vertices.

Recall that for a finite group $X$, any indecomposable $FX$-module $Q$
has a vertex. This is, by definition, a subgroup $Y$ of $X$ such that
$Q$ is $Y$-projective and which is minimal with this property.  (A
module is $Y$-projective if it is a summand of $S\tensor_{FY}FX$ for
some $FY$-module $S$).  Any two vertices are conjugate in $X$, and
moreover they are $p$-subgroups \cite[Proposition~3.10.2]{benson98}.
The modules with vertex $\{1\}$ are precisely the projective modules.

By Proposition~\ref{p-permu}, $S^p(U)$ is $\Sym_p$-projective.  A
Sylow $p$-subgroup $D$ of $\Sym_p$ is cyclic of order $p$.  Hence an
indecomposable summand of $S^p(U)$ is non-projective if and only if it
has vertex $D$.

\begin{lemma} \label{step1}
  There is a {1-1} correspondence, multiplicity-preserving, between the
  indecomposable summands of $S^p(U)$ with vertex $D$ and the
  indecomposable summands of $M$ with vertex $D$ where
  $$
  M=F\tensor_{F\Sym_p} e^{\# p}FH.
  $$
\end{lemma}

\begin{proof}
  Recall that $S^p(U)\cong M\tensor_{FH} F\Sym_n$ where
  $H=\Sym_k\wr\Sym_p$.  To prove the statement it is sufficient to
  establish the following.
  
  If $X$ is an indecomposable summand of $M$ with vertex $D$ then
  $X\tensor_{FH}F\Sym_n$ has a unique indecomposable summand
  $\tilde{X}$ with vertex $D$; and moreover $X\mapsto \tilde{X}$ is a {1-1}
  correspondence.
  
  Let $N_1 = N_H(D)$ and $N=N_{\Sym_n}(D)$.  The Green
  correspondence~\cite[Theorem~3.12.2]{benson98} provides us with a
  {1-1} correspondence between the indecomposable $FH$-modules (or
  $F\Sym_n$-modules) with vertex $D$ and the indecomposable
  $FN_1$-modules (or $FN$-modules) with vertex $D$. So we are done if
  we can show that if $Q$ is the Green correspondent of $X$ in $N_1$
  then $Q\tensor_{FN_1}FN$ is indecomposable, and that non-isomorphic
  Green correspondents induce to non-isomorphic $FN$-modules. (The
  module $Q\tensor_{FN_1}FN$ has then automatically vertex $D$.)
  
  Since $X$ is a $p$-permutation module with vertex $D$, one knows
  that $D$ acts trivially on $Q$ and moreover $Q$ is indecomposable
  projective as a module for $N_1/D$. Then $Q\tensor_{FN_1}FN$ is
  still trivial as a module for $D$, and it is projective as a module
  for $N/D$. To complete the proof we exploit the general fact that
  indecomposable projective modules are in 1-1 correspondence with
  their simple quotients.
  
  Analysing the groups $N_1$ and $N$ we will show below (in the
  following Lemma) that
  $$
  N = N_1B
  $$
  where $B$ is a $p$-group which is normal in $N$.
  
  Recall that a normal $p$-subgroup acts trivially on all simple
  modules, therefore we have a {1-1} correspondence between simple
  modules of $N$ and simple modules of $N_1$, by restriction.  By
  Frobenius reciprocity, it follows that
  $$
  \Hom_N(Q\tensor_{FN_1}FN, L)\cong \Hom_{N_1}(Q,L)
  $$
  for any simple $N$-module $L$ and this is $F$ precisely if $L$ is
  the simple quotient of $Q$, and zero otherwise.
\end{proof}

\begin{lemma}
  Let $N_1 = N_H(D)$ and $N=N_{\Sym_n}(D)$. Then $N=N_1B$ where $B$ is
  a $p$-group which is normal in $N$.  Moreover $N_1/D$ is isomorphic
  to the direct product of $\Sym_k$ with a cyclic group of order
  $p-1$.
\end{lemma}

\begin{proof} 
  Observe first that the centraliser $C_H(D)$ of $D$ in $H$ is
  isomorphic to a direct product $\Sym_k\times D$. We view this in two
  ways,
  
  (i) as a subgroup of $\Sym_k\wr \Sym_p$ where the first factor of
  the above direct product is the diagonal $\Delta(\Sym_k)$ in the
  base group; and
  
  (ii) as a subgroup of $N_{\Sym_p}(C_p)\wr \Sym_k$ where the second
  factor of the direct product above is contained in the diagonal of
  the base group.
  
  We get from (i) that $N_1 \cong \Delta(\Sym_k)\times N_{\Sym_p}(D)$.
  Moreover, we get from (ii) that $C_{\Sym_n}(D)$ is the semi-direct
  product of $B$ with $\Sym_k$ where $B$ is the base group of $C_p\wr
  \Sym_k$, and that $N$ is generated by $C_{\Sym_n}(D)$ and
  $N_{\Sym_p}(D)$ (contained in the base group) which acts diagonally.
  Hence $B$ is normal in $N$ and $N= BN_1$.
  
  The description of $N_1$ implies directly the statement on $N_1/D$.
\end{proof}

We will apply Brou{\'e}'s correspondence theorem for $p$-permutation
modules which is described now; here $X$ can be an arbitrary finite
group again. Assume $W$ is a $p$-permutation module of $X$ and $P$ is
a $p$-subgroup of $X$.  Set
$$
W(P):= W^P/\sum_{Q<P} \textrm{Tr}_Q^P(W^Q),
$$
where $W^R$ denotes the space of fixed points in $W$ of any
subgroup $R$ of $X$ and where $\textrm{Tr}_Q^P$ is defined as
$\textrm{Tr}_Q^P(m) =\sum_i m g_i$, the sum taken over a transversal
of $Q$ in $P$.  Then $W(P)$ is a module for $N_G(P)$ on which $P$ acts
trivially, and hence is a module for the factor group $N_G(P)/P$.  As
a vector space, $W(P)$ is isomorphic to the span of the fixed points
of $P$ in a given permutation basis.

In~\cite{broue85}, Brou{\'e} proved that there is a {1-1}
correspondence, multiplicity-preser\-ving, between the indecomposable
summands of $W$ with vertex $P$ and the indecomposable summands of
$W(P)$ which are projective as $N_G(P)/P$-modules.
 
Applied to the $p$-permutation module $M$ of $H$ and combined with
Lemma~\ref{step1}, this result implies that there is a {1-1}
correspondence between non-projective indecomposable summands of
$S^p(U)$ and the indecomposable summands of $M(D)$ which are
projective as $N_H(D)/D$-modules.  This group is isomorphic to the
direct product of $\Sym_k$ with a cyclic group of order $p-1$ whose
group algebra over $F$ is semi-simple with $1$-dimensional simple
modules.  The proof of Theorem~\ref{param} can thus be completed using
the following result.

\begin{lemma}
  The $N_1/D$-module $M(D)$ is isomorphic to $U$ as an
  $\Sym_k$-module. Moreover, the cyclic group of order $p-1$ acts
  trivially on $M(D)$.
\end{lemma}

\begin{proof}
  Let $B$ and $B'$ be bases of $U=eF\Sym_k$ and $(1-e)F\Sym_k$,
  respectively, so that $B\cup B'$ is a basis of $F\Sym_k$.  Then the
  induced module $F\tensor_{F\Sym_p}FH$ has basis
  $
  \setofall{s_p\tensor (v_1\# \ldots \# v_p)}{v_1,\ldots,v_p\in B\cup  B'}$.  
  It follows that
  $$
  \setofall{s_p\tensor (b_1\# \ldots \# b_p)}{b_1,\ldots,b_p\in B}
  $$
  is a basis of $M$ where $M=F\tensor_{F\Sym_p}e^{\# p}FH$.  This is a
  permutation basis under the action of $\Sym_p$.
  
  If, in particular, $\zeta\in D$ is a $p$-cycle which in its action
  on $kp$ points permutes the supports of the factors cyclically, then
  $$
  \Big(s_p\tensor (b_1\#b_2\# \ldots \# b_p)\Big)\zeta 
  = 
  s_p\tensor (b_p\# b_1\# \ldots \# b_{p-1})
  $$
  for all $b_1,\ldots,b_p\in B$.  We deduce that a basis vector of
  $M$ is fixed by $D$ if and only if it is of the form $s_p\tensor
  (b\# b\# \ldots \# b)$ for some $b\in B$.  Such element is also
  fixed under the cyclic group of order $p-1$ normalising $D$ in
  $\Sym_p$ (which proves the last part of the Lemma).
  
  As a consequence, there is an obvious vector space isomorphism
  $\psi: U \to M(D)$ taking $b$ to the coset of $s_p\tensor (b\#
  b\#\ldots \# b)$ for all $b\in B$.  We want to compare the actions
  of $\Sym_k$.  On the one hand, we have for $\pi\in \Sym_k$
  $$
  \psi(b\pi) 
  = 
  \psi \Big(\sum_{b'\in B} c(b,b')\,b'\Big) 
  = 
  \sum_{b'\in B} c(b,b')\,\psi(b')
  $$
  with coefficients $c(b,b')\in F$. On the other hand,
  $$
  \psi(b)\pi 
  = 
  s_p\tensor (b\pi\# b\pi\# \ldots \# b\pi)
  = 
  \sum_{b'\in B} c(b,b')^p\;s_p\tensor (b'\# \ldots \# b') + (*) 
  $$
  where $(*)$ belongs to the span of orbit sums with orbits of size
  $>1$ and which is zero in $M(D)$.
  
  So $M(D)$ is isomorphic to the Frobenius twist of $U$, obtained by
  composing the corresponding matrix representation with the map
  $$
  \Big(c(b,b')\Big) \mapsto \Big(c(b,b')^p\Big).
  $$
  But since all projective modules for symmetric groups are defined
  over the prime field, for any projective module $P$, its Frobenius
  twist is isomorphic to $P$.  So the Lemma is proved.
\end{proof}

Combining Theorem~\ref{param} with Theorem~\ref{kp} and applying the
Heller operator, yields Theorem~\ref{sym-haupt} mentioned in the
introduction.

Recall that an $F\Sym_n$-module $M$ is \emph{liftable} if there exists
an indecomposable $\mathcal{O}\Sym_n$-module $\tilde{M}$ whose
reduction modulo $p$ is $M$; here $\mathcal{O}$ is a complete discrete
valuation ring with residue field $F$. The above considerations also
have the following consequence.

\begin{corollary} \label{liftable}
  Let $k$ be a positive integer not divisible by $p$, then any
  non-projective indecomposable summand of $L_{pk}$ is liftable and
  has an associated complex character.
\end{corollary}

\begin{proof}
  Any $p$-permutation module is liftable, by Scott's
  theorem~\cite[3.11.3]{benson98}, hence every indecomposable summand
  of $S^p(L_k)$ is liftable.  Furthermore, an $\Omega$-translate of a
  liftable module is liftable, by a standard argument.  The claim
  follows from Theorem~\ref{kp} and Proposition~\ref{p-permu}.
\end{proof}

In concluding this section, we give a more detailed account on the Lie
module $L_{2p}$, for any odd prime $p$.  This is a case where the
modular representation theory of the symmetric group is sufficiently
well understood, and we have a complete description of $S^p(L_2)$.

\begin{theorem} \label{2p}
  Let $p$ be an odd prime, then $L_{2p}$ has a unique non-projective
  indecomposable summand. It belongs to the principal block
  and has character
  \begin{eqnarray*}
    \lefteqn{%
      \chi^{(2,1^{2p-2})}
      + 
      \chi^{(3,2^{p-2},1)}}\\[1mm]
    & + &
    \sum_{i=2}^{(p-1)/2}
    \Big(
    \chi^{(2i+1,2i,2^{p-2i-1},1)}
    + 
    \chi^{((2i-1)^2,2^{p-2i+1})}
    + 
    \chi^{(2i,2i-1,2^{p-2i},1)}
    \Big),
  \end{eqnarray*}
  where $\chi^\mu$ denotes the irreducible character of $\Sym_{2p}$
  corresponding to $\mu$, for any partition $\mu$ of $2p$.
\end{theorem}

Since the character of $L_{2p}$ is known, one knows the character of
the projective part of $L_{2p}$ as well. Since a projective module is
determined by its character one can deduce the complete direct sum
decomposition of $L_{2p}$, at least in principle.

The proof of Theorem~\ref{2p} builds on the following result on
$S^p(L_2)$.

\begin{proposition} \label{spl2}
  Let $p$ be an odd prime, then $S^p(L_2)$ has a unique non-projective
  indecomposable summand $Q$, namely its principal block component.
  The character of $Q$ is 
  $$
  \chi^{(1^{2p})}
  +
  \sum_{i=1}^{(p-1)/2} \chi^{((2i+1)^2,2^{p-2i+1})}.
  $$ 
\end{proposition}

\begin{proof}
  Let $H=\Sym_2\wr \Sym_p$ and
  $M=F\tensor_{F\Sym_p}\omega^{(2,\ldots,2)}FH$.  We denote the sign
  representation of $F\Sym_{2p}$ by $\sgn$, then $M$ is isomorphic to
  the sign representation $\sgn|_{FH}$ of $FH$, since the idempotent
  $\frac{1}{2}\omega_2$ generates the sign representation of $\Sym_2$.
  Applying Proposition~\ref{p-permu} and the tensor identity, it
  follows that
  $$
  S^p(L_2)\tensor\sgn
  \cong
  (M\tensor\sgn|_{FH})\tensor_{FH} F\Sym_n
  =
  F\tensor_{FH} F\Sym_n\,.
  $$
  By definition, $H$ is the centraliser of a fixed point free
  involution, $\tau$ say, hence the module $S^p(L_2)\tensor \sgn$ is
  isomorphic to the permutation module of the symmetric group
  $\Sym_{2p}$ on the conjugacy class of $\tau$.
  
  This module was studied in detail by
  Wildon~\cite[Chapter~5]{wildon04}. The component in the principal
  block is indecomposable, and this is preserved by tensoring with the
  sign representation.  However, by Theorem~\ref{param}, $S^p(L_2)$
  has a unique non-projective indecomposable summand $Q$.  Its Green
  correspondent belongs to the principal block, hence so does $Q$.
  This proves the first claim.
  
  The character of $S^p(L_2)$ is $ \sum\chi^\mu $, where the sum is
  taken over all $\mu$ such that the multiplicity of each part of
  $\mu$ is even (see, for example, \cite[Theorem~4.1.1]{wildon04}).
  In order to determine its principal block component, it is
  convenient to use abacus notation for characters (see
  \cite[p.~78]{martin93}, or \cite{scopes91}). Any partition $\lambda$
  of $2p$ can be represented on an abacus with $p$ runners and with
  $2p$ beads. Then $\lambda$ lies in the principal block if and only
  if the abacus display has two gaps (counting the gaps on each runner
  which are above the last bead). We write $\abacus{i}{j}$ if the gaps
  are on runner(s) $i, j$ and write $\langle\,i\,\rangle$ if there is
  a gap of size $2$ on runner $i$.  These give a complete list of
  partitions whose character belongs to the principal block.
  
  For example, if $p=5$, then $\abacus{1}{3}$ denotes the abacus
  display
  $$
  \left(
    \begin{array}{ccccc}
      \bullet & \bullet & \bullet & \bullet & \bullet \\[1mm]
      \cdot   & \bullet & \cdot   & \bullet & \bullet \\[1mm]
      \bullet & \cdot   & \bullet & \cdot   & \cdot   \\[1mm]
    \end{array}
  \right).
  $$
  This represents the partition $\lambda=(3,2,2,2,1)$ (count the
  number of gaps before each bead, reading row-wise from top to bottom).
  
  By a straightforward case-by-case analysis one now shows that the
  character of $Q$ is
  \begin{equation}  \label{aba-char}
    \abacus{1}{1} + \sum_{i=1}^{(p-1)/2} \abacus{2i}{2i+1}
  \end{equation}
  and this translates directly into the statement.
\end{proof}

\begin{proof}[Proof of Theorem~\ref{2p}]
  By Theorem~\ref{kp} and Proposition~\ref{spl2}, the Lie module
  $L_{2p}$ has a unique non-projective indecomposable summand which
  belongs to the principal block, namely the Heller translate of $Q$.
  
  We use abacus notation for characters and irreducible modules in the
  principal block component again. From~\cite{wildon04}, $Q$ has top
  composition factors
  $$
  \bigoplus_{i=1}^{p-1/2} D(\abacus{ 2i+1, 2i }).
  $$
  The decomposition matrix of the block is known~\cite{martin93}.
  The projective
  cover of $D(\abacus{ 2i+1, 2i })$ has character
  $$
  \abacus{2i+1}{2i} 
  + 
  \abacus{2i+1}{2i-1} 
  + 
  \abacus{2i-1}{2i-2} 
  + 
  \abacus{2i}{2i-2}
  $$
  if $i\ge 2$, and
  $\abacus{3}{2}+\abacus{3}{1}+\abacus{2}{2}+\abacus{1}{1}$ if $i=1$.
  Comparison with the character of $Q$ as given in~\eqref{aba-char},
  shows that the character of the non-projective indecomposable
  summand of $L_{2p}$ is
  $$
  \abacus{2}{2} 
  + 
  \abacus{3}{1} 
  + 
  \sum_{i=2}^{(p-1)/2}
  \Big(
  \abacus{2i+1}{2i-1} 
  + 
  \abacus{2i-1}{2i-2} 
  + 
  \abacus{2i}{2i-2}
  \Big),
  $$
  which readily translates into the claim.
\end{proof}

For general $n=kp$ with $k>2$ not divisible by $p$, one does not know
the representation theory of symmetric group and the character of
$S^p(L_k)$ in sufficient detail in order to derive a result like
Theorem~\ref{2p} in this way.

\section{Lie powers in degrees not divisible by $p^2$}  \label{6}

Let $F$ be an infinite field of prime characteristic $p$ and let $V$
be an $r$-dimensional vector space over $F$.  In this section, we
consider Lie powers $L_{pk}(V)$ of $V$ where $k$ is not divisible by
$p$.

If $n$ is an arbitrary positive integer, then $V^{\tensor n}$ is an
$(\Sym_n,\GL(V))$ bimodule.  The action of $\GL(V)$ from the right is
the diagonal action, while for the $\Sym_n$-action from the left, we
have
$$
\pi
(v_1\tensor\cdots\tensor v_n)
=
v_{1\pi}\tensor\cdots\tensor v_{n\pi}
$$
for all $\pi\in\Sym_n$ and $v_1,\ldots,v_n\in V$.  

For $r\ge n$, the Schur algebra $S(r,n)$ contains an idempotent $\xi$
such that the algebra $\xi S(r,n)\xi$ is canonically isomorphic to the
group algebra $F\Sym_n$ and we identify these algebras
\cite[(6.1d)]{green80}\footnote{%
  Note that the roles of $r$ and $n$ here and left and right action
  are exchanged in comparison with~\cite{green80}.}.  The Schur functor
(denoted by $f$ in \cite{green80}) takes an $S(r,n)$-module $M$ to the
$F\Sym_n$-module $M\xi$. The tensor space $V^{\tensor n}$ is
isomorphic to $\xi S(r,n)$ as an $(F\Sym_n,S(r,n))$
bimodule~\cite[(6.4f)]{green80}. Hence the left adjoint of the Schur
functor (which we denote by $g_\tensor$) takes a right
$F\Sym_n$-module $U$ to $U\tensor_{F\Sym_n} V^{\tensor n}$.

The functor $g_{\tensor}$ is right exact but not exact.  However, it
is possible to lift the short exact sequence of $\Sym_n$-modules given
in Theorem~\ref{kp} to a short exact sequence of $\GL(V)$-modules and
to parametrise the indecomposable summands of $L_{pk}(V)$ which are
not summands of $V^{\tensor n}$ accordingly.

We will use the following tool, to translate between the symmetric
group and the general linear group.

\begin{lemma} \label{lifting}
  Let $n$ be a positive integer and assume that
  $$
  0
  \lra
  fF\Sym_n
  \stackrel{\alpha}{\lra} 
  e F\Sym_n
  \lra
  U
  \lra 
  0
  $$
  is a short exact sequence of $\Sym_n$-right modules, where
  $e,f\in F\Sym_n$ satisfy $ef=f$ and $e^2=e$, and where $\alpha$ is
  left multiplication with $e$.  Then there is a short exact sequence
  of $\GL(V)$-modules
  $$
  0
  \lra
  fV^{\tensor n}
  \lra
  e\tensor_{F\Sym_n} V^{\tensor n}
  \lra
  U\tensor_{F\Sym_n} V^{\tensor n}
  \lra 
  0.
  $$
\end{lemma}

\begin{proof}
  There is an embedding of $\GL(V)$-modules $\alpha':fV^{\tensor n}\to
  eV^{\tensor n} $ provided by left action of $e$.  Furthermore,
  application of $g_\tensor$ to the short exact sequence of
  $\Sym_n$-modules gives the exact sequence
  $$
  f\tensor_{F\Sym_n} V^{\tensor n}
  \stackrel{\gamma}{\lra}
  e\tensor_{F\Sym_n} V^{\tensor n}
  \lra
  U\tensor_{F\Sym_n} V^{\tensor n}
  \lra
  0
  $$
  of $\GL(V)$-modules.  The multiplication map
  $F{\Sym_n}\tensor_{F\Sym_n}V^{\tensor n} \to F{\Sym_n}V^{\tensor n}$
  restricts to epimorphisms $\beta:f\tensor_{F\Sym_n} V^{\tensor n}\to
  fV^{\tensor n}$ and $\delta:e\tensor_{F\Sym_n} V^{\tensor n}\to
  eV^{\tensor n}$ and gives rise to a commutative diagram
  $$
  \xymatrix@C=5ex@R=6ex{%
    &
    f\tensor_{F\Sym_n} V^{\tensor n}
    \ar@{->}[r]^-{\mbox{$\gamma$}}  
    \ar@{->}[d]_-{\mbox{$\beta$}} \;
    &
    e\tensor_{F\Sym_n} V^{\tensor n}
    \ar@{->}[d]_-{\mbox{$\delta$}}\;\;
    \ar@{->}[r]
    &
    U\tensor_{F\Sym_n} V^{\tensor n}
    \ar@{->}[r]
    &
    0\\ 
    0
    \ar@{->}[r]
    &
    fV^{\tensor n}
    \ar@{->}[r]_-{\mbox{$\alpha'$}}  
    & 
    eV^{\tensor n}
    &
    &
  }
  $$
  But $\delta$ is an isomorphism, since $e$ is an idempotent.  It
  follows that $ \ker\,\gamma=\ker\,\beta $ and that
  $\delta^{-1}\alpha'$ is an embedding of $fV^{\tensor n}$ into
  $e\tensor_{F\Sym_n} V^{\tensor n}$ such that
  $\delta^{-1}\alpha'\beta=\gamma$.
\end{proof}

\begin{corollary} \label{GL-kp}
  If $n=pk$ such that $k$ is not divisible by $p$, then there is a
  short exact sequence of $\GL(V)$-modules
  $$
  0
  \lra
  L_n(V)
  \lra
  e_n\tensor_{F\Sym_n} V^{\tensor n}
  \lra
  S^p(L_k(V))
  \lra 
  0,
  $$
  where $S^p(L_k(V))$ denotes the $p$-th symmetric power of
  $L_k(V)$.
\end{corollary}

\begin{proof}
  This follows from Theorem~\ref{kp} and Lemma~\ref{lifting}, applied
  to $f=\omega_n$, $e=e_n$ and $U=S^p(L_k)$, we just have to identify
  the cokernel.  Let $H=\Sym_k\wr\Sym_p$, then $\omega^{(k,\ldots,k)}$
  is (up to the factor $1/k^p$) an idempotent in $FH$.  Hence, by
  Proposition~\ref{p-permu},
  \begin{eqnarray*}
    S^p(L_k)\tensor_{F\Sym_n} V^{\tensor n}
    & \cong &
    \Big((F\tensor_{F\Sym_p}\omega^{(k,\ldots,k)}FH)\tensor_{FH} F\Sym_n\Big)
    \tensor_{F\Sym_n} 
    V^{\tensor n}\\[1mm]
    & \cong &
    (F\tensor_{F\Sym_p} \omega^{(k,\ldots,k)} F\Sym_n)
    \tensor_{F\Sym_n} 
    V^{\tensor n}\\[1mm]
    & \cong &
    F\tensor_{F\Sym_p} \omega^{(k,\ldots,k)} V^{\tensor n}\\[1mm]
    & \cong &
    F\tensor_{F\Sym_p} L_k(V)^{\tensor p} \\[1mm]
    & \cong &
    S^p(L_k(V))
  \end{eqnarray*}
  as desired.
\end{proof}

Note that, as the proof shows, the image under $g_\tensor$ of
$S^p(L_k)$ is isomorphic to $S^p(L_k(V))$. Therefore the image under
the Schur functor $f$ of $S^p(L_k(V))$ is isomorphic to $S^p(L_k)$,
since $f\circ g_\tensor$ is naturally equivalent to the identity.

To prove Theorem~\ref{GL-haupt}, we will use the following general
observations on the functor $g_{\tensor}$.

\begin{proposition}
  Let $n$ be a positive integer such that $r\ge n$.  Then
\begin{enumerate}
  \item If $M$ is an indecomposable $F\Sym_n$-module, then
    $M\tensor_{F\Sym_n}V^{\tensor n}$ is indecomposable.
  \item If $M_1$ and $M_2$ are indecomposable $F\Sym_n$-modules which
    are not isomorphic, then $M_1\tensor_{F\Sym_n}V^{\tensor n}$ and
    $M_2\tensor_{F\Sym_n}V^{\tensor n}$ are not isomorphic.
  \end{enumerate}
\end{proposition}

\begin{proof}
  Let $S=S(r,n)$ be the Schur algebra and let $f$ be the Schur
  functor. Recall that $g_{\tensor}$ is left adjoint to $f$ and that
  this functor followed by the Schur functor is naturally equivalent
  to the identity.  So
  $$
  \Hom_S(g_{\tensor} M, g_{\tensor} M) 
  \cong 
  \Hom_{F{\cal S}_n}(M, f\circ g_{\tensor} M) 
  \cong 
  \Hom_{F{\cal S}_n}(M, M),
  $$
  by adjointness, even an isomorphism of algebras. It follows that
  $M$ is indecomposable if and only if $g_{\tensor} M$ is
  indecomposable.
  
  To prove the second part, let $M_1$ and $M_2$ so that
  $g_{\tensor} M_1$ and $g_{\tensor} M_2$ are isomorphic, then
  also
  $
  M_1 \cong f\circ g_{\tensor} M_1 \cong f\circ g_{\tensor} M_2 \cong M_2
  $.
\end{proof}

For $r\ge n$, any summand of the tensor space $V^{\tensor n}$ is
projective and injective as a module for the Schur algebra $S(r,n)$
(see, for example,~\cite[p.~94]{donkin98}). As a consequence, the
Heller operator $\Omega$ gives a {1-1} correspondence between
indecomposable non-projective quotients of $V^{\tensor n}$ and
indecomposable non-injective submodules of $V^{\tensor n}$.

Furthermore, if $M$ is a quotient of $V^{\tensor n}$ then $M$ is not
projective if and only if it is not projective and injective.  One
direction is clear. For the converse, assume $M$ is projective then it
is a summand of $V^{\tensor n}$ and hence is also injective.
Similarly, a submodule of $V^{\tensor n}$ is not injective if and only
if it is not projective and injective.

\begin{proposition} \label{GL-param}
  Let $n=pk$ such that $k$ is not divisible by $p$ and assume $r\ge
  n$.
  \begin{enumerate}
  \item The functor $g_{\tensor}$ gives a {1-1} correspondence,
    multiplicity-preserving, between the non-projective indecomposable
    summands of the $F\Sym_n$-mo\-dule $S^p(L_k)$ and the
    indecomposable summands of the $S(r,n)$-module $S^p(L_k(V))$ which
    are not projective and injective.
  
  \item The Heller operator $\Omega$ gives a {1-1} correspondence,
    multiplicity-preserv\-ing, between the indecomposable summands of
    $S^p(L_k(V))$ which are not projective and injective and the
    indecomposable summands of $L_n(V)$ which are not projective and
    injective, both as modules for $S(r,n)$.
  \end{enumerate}
\end{proposition}

\begin{proof}
  1. By the preceding proposition we get the {1-1} correspondence
  between indecomposable summands of $S^p(L_k)$ and $S^p(L_k(V))$.
  Moreover, an indecomposable $F{\cal S}_n$-module $M$ is projective
  if and only if $g_{\tensor} M$ is an indecomposable summand of
  $V^{\tensor n}$ since $g_\tensor$ takes $F\Sym_n$ to $V^{\tensor
    n}$. This completes the proof of the first part.

  2. Consider the short exact sequence
  $$
  0 \to L_n(V) \to e_nF\Sym_n \tensor V^{\tensor n} \to S^p(L_k(V))\to 0.
  $$
  The middle is isomorphic to a summand of $V^{\tensor n}$.
  As explained above, $\Omega$ induces a {1-1} correspondence as stated
  in the second part.
\end{proof}

It was shown in~\cite{donkin-erdmann98} that there is a {1-1}
correspondence between the indecomposable summands of $L_k$ and those
of $L_k(V)$.  Explicitly, for $k$ not divisible by $p$, the module
$L_k(V)$ is a summand of $V^{\tensor k}$. The indecomposable summands
of $V^{\tensor k}$ are parametrised by $p$-regular partitions
$\lambda$ of $k$, as $T(\lambda)$, where $T(\lambda)$ has highest weight $\lambda$ and is also known as tilting module.
Via the Schur functor, $T(\lambda)$
corresponds to the projective $\Sym_k$-module $P(\lambda)$ with simple
quotient $D^\lambda$ labelled by  the partition $\lambda$.

Theorem~\ref{GL-haupt} is now an immediate consequence
of Proposition~\ref{GL-param} and Theorem~\ref{param}, applied to
$U=L_k$.

\begin{remark}
  We have assumed that the dimension $r$ of $V$ should be $\ge pk$.
  One might ask what can be deduced by truncating the exact sequence
  $$
  0
  \to 
  L_{pk}(V) 
  \to 
  e_n
  \tensor_{F\Sym_n} V^{\tensor n} 
  \to 
  S^p(L_k(V)) 
  \to 
  0.
  $$
  See \cite[Section~6.5]{green80} for details on truncation.  Take
  $d < r$ and take a subspace $E$ of $V$ of dimension $d$ (with basis
  a subset of the canonical basis of $V$). There is an idempotent
  $\zeta \in S(r,n)$ such that $\zeta S(r,n)\zeta \cong S(d,n)$. The
  functor $(-)\zeta$ is exact. It takes $V^{\tensor n} $ to
  $E^{\tensor n}$ and $L_{pk}(V)$ to $L_{pk}(E)$. So we get an exact
  sequence
  $$
  0
  \to 
  L_{pk}(E)
  \to 
  e_n\tensor_{F\Sym_n} E^{\tensor n}
  \to 
  S^p(L_k(V))\zeta
  \to 
  0.
  $$
  However, $E^{\tensor n}$ is not projective and injective in
  general as a module for $S(d,n)$.  Moreover, the kernel of this
  sequence can sometimes be isomorphic to an indecomposable summand of
  $E^{\tensor n}$.
  
  For example, let $p=3 =n$ and $r=2$. As mentioned in
  \cite[Example~3.6]{donkin-erdmann98}, the module $L_3(E)$ is
  isomorphic to $L_2(E)\tensor E$ and is therefore isomorphic to a
  summand of $E^{\tensor 3}$.  (In fact, one can show that this
  summand is neither projective nor injective as a module for the
  Schur algebra $S(2,3)$.)
\end{remark}

In concluding this section, we recover a result of Bryant and
St{\"o}hr on the $p$th Lie power.

Let $L(V)''$ denote the second derived algebra of $L(V)$, then the
quotient $M(V)=L(V)/L(V)''$ is a free metabelian Lie algebra.  It was
shown in~\cite{bryant-stoehr00} that
$$
L_p(V) \cong B_p(V) \oplus M_p(V),
$$
where $B_p(V)=L(V)''\cap V^{\tensor p}$ and
$M_p(V)=(L_p(V)+L(V)'')/L(V)''$.

\begin{corollary}[{\cite[Theorem~3.1]{bryant-stoehr-pre}}] \label{GL-p}
  $B_p(V)$ is a summand of $V^{\tensor p}$.
\end{corollary}

\begin{proof}
  Assume first that $r\ge p$ and consider the short exact sequence
  $$
  0
  \to
  M_p(V)
  \to
  V\tensor S^{p-1}(V)
  \to
  S^p(V)
  \to 
  0.
  $$
  (see, for instance,~\cite{hanne-sto90}) This is non-split since,
  for example, the Schur functor takes the middle to the natural
  permutation module of ${\cal S}_p$ which is indecomposable.
  Therefore $M_p(V)$ is not injective.
  
  By Theorem~\ref{GL-haupt}, we know that $L_p(V)$ has a unique
  indecomposable summand which is not projective and injective.  It
  follows that $B_p(V)$ is injective, hence is a summand of
  $V^{\tensor p}$.
  
  For $r<p$, apply the truncation method as described in the previous
  remark.
\end{proof}

\section{Factorisation of Lie resolvents} \label{7}

Let $G$ be a group and $F$ be a field. The Green ring $R_{FG}$ of $G$
over $F$ has basis the isomorphism classes of finite-dimensional
indecomposable $FG$-modules with addition and multiplication arising
from direct sums and tensor products. If $V$ is a finite-dimensional
$FG$-module, we also write $V$ for the corresponding element in
$R_{FG}$. So, for instance, $V^n\in R_{FG}$ is the isomorphism class
of the $n$th tensor power of $V$.

If $V$ is a finite-dimensional $FG$-module, then the $n$th Lie power
$L_n(V)$ may be regarded as a module for $FG$ in a natural way.
Recently, Bryant \cite{bryant-2,bryant-adams,bryant-1} introduced the
\emph{Lie resolvents} $\phi_{FG}^n:R_{FG}\to R_{FG}$, $n\ge 1$, to
study the structure of $L_n(V)$.  These can be described by
$$
\phi^n_{FG}(V)
=
\sum_{d\divides n} \mu(n/d) dL_d(V^{n/d})
$$
for all $n$ and $V$, where $\mu$ denotes the M{\"o}bius function.
By M{\"o}bius inversion, this is equivalent to
$$
L_n(V)
=
\frac{1}{n}\,\sum_{d\divides n} \phi_{FG}^d(V^{n/d})
$$
for all $n$ and all modules $V$.  Thus complete knowledge of the
Lie resolvents $\phi^n_{FG}$ yields a description of $L_n(V)$ for each
$V$, up to isomorphism.

Most strikingly, the Lie resolvents are \emph{linear endomorphisms} of
$R_{FG}$ (see \cite[Corollary~3.3]{bryant-adams}).  Let $p$ denote the
characteristic of $F$ (which may or may not be zero), then $
\phi_{FG}^{k'k} = \phi_{FG}^{k'} \circ \phi_{FG}^k $ for all coprime
positive integers $k$, $k'$ not divisible by $p$ (see Theorem~5.4 and
Corollary~6.2 in \cite{bryant-adams}).  If the characteristic $p$ of
$F$ is positive and $G$ is finite with Sylow $p$-subgroups of order
$p$, then there is also the identity
$$
\phi_{FG}^{p^mk}
=
\phi_{FG}^{p^m}
\circ
\phi_{FG}^k
$$
for all positive integers $m$, $k$ such that $k$ is not divisible
by $p$ (see \cite[Corollary~1.2]{bryant-2}). The question was raised
in \cite{bryant-2} whether such factorisation rule might hold for
arbitrary groups $G$.

We shall give here an answer in case $m=1$.

\begin{theorem} \label{pk-gl-pk}
  Let $G$ be a group and $F$ be a field of prime characteristic $p$.
  Then, for any positive integer $k$ not divisible by $p$,
  $$
  \phi_{FG}^{pk}
  =
  \phi_{FG}^p\circ\phi_{FG}^k.
  $$  
\end{theorem}

\begin{proof}
  Let $n=pk$. It suffices to prove 
  \begin{equation} \label{fac-mod}
    \phi_{FG}^{pk}(V) = \phi_{FG}^p(\phi_{FG}^k(V))
  \end{equation}
  for any finite-dimensional $FG$-module $V$, by linearity.  In fact,
  it suffices to consider the case where $F$ is infinite, $G=\GL(V)$
  and $V$ has dimension $r\ge n$, by~\cite[Lemma~2.4]{bryant-2} and
  arguments completely analogous to those given in the proof
  of~\cite[Theorem~5.4]{bryant-adams}.  We write $\phi$ for
  $\phi_{FG}$.
  
  The equality~\eqref{fac-mod} holds if $F$ has characteristic zero,
  by Corollary~6.2 and Theorem~5.4 in \cite{bryant-adams}.  In
  particular, $\phi^{pk}(V)$ and $\phi^p(\phi^k(V))$ have the same
  formal character.
  
  We shall now show~\eqref{fac-mod} for all $V$ of dimension $r\ge n$,
  by induction on $k$.  For $k=1$, this follows from
  $\phi^1(V)=L_1(V)=V$.  Let $k>1$, then inductively
  \begin{eqnarray*}
    nL_n(V)
    & = &
    \sum_{d\divides k} (\phi^{pd}(V^{k/d})+\phi^d(V^{pk/d}))\\[1mm]
    & = &
    \phi^{pk}(V)+\phi^k(V^p)
    +
    \sum_{d\divides k,\,d\neq k} (\phi^{pd}(V^{k/d})+\phi^d(V^{pk/d}))\\[1mm]
    & = &
    \phi^{pk}(V)-\phi^p(\phi^k(V))
    +
    \phi^p\Big(\sum_{d\divides k} \phi^d(V^{k/d})\Big)
    +
    \sum_{d\divides k} \phi^d(V^{pk/d})\\[1mm]
    & = &
    \phi^{pk}(V)-\phi^p(\phi^k(V))
    +
    nL_p(L_k(V))-kL_k(V)^p
    +
    kL_k(V^p).
  \end{eqnarray*}
  As a consequence, $\phi^{pk}(V)-\phi^p(\phi^k(V)) $ is a linear
  combination of $L_n(V)-L_p(L_k(V))$, $L_k(V^p)$ and $L_k(V)^p$ in
  $R_{FG}$. The last two are summands of $V^{\tensor n}$
  since $p$ does not divide $k$.
  
  Concerning the first one, note that $L_p(L_k(V))$ embeds naturally
  into $L_n(V)$.  By Theorem~\ref{GL-kp}, there are the two short
  exact sequences of $\GL(V)$-modules
  $$
  0
  \lra
  L_p(L_k(V))
  \lra
  T_1
  \lra
  S^p(L_k(V))
  \lra 
  0
  $$
  and
  $$
  0
  \lra
  L_n(V)
  \lra
  T_2
  \lra
  S^p(L_k(V))
  \lra 
  0,
  $$
  where $T_1=e_p\tensor_{F\Sym_p} (L_k(V))^{\tensor p}$ and
  $T_2=e_n\tensor_{F\Sym_n} V^{\tensor n}$.  The modules $T_1$ and
  $T_2$ are summands of $V^{\tensor n}$, thus they are projective as
  modules for the Schur algebra $S(r,n)$.  Applying Schanuel's Lemma,
  we have $L_p(L_k(V))+T_2=T_1+L_n(V)$ in $R_{FG}$.
  
  This shows that $\phi^{pk}(V)-\phi^p(\phi^k(V))$ is an integer linear
  combination of summands of $V^{\tensor n}$, which allows us to
  deduce~\eqref{fac-mod} from the equality of the corresponding
  characters \cite{donkin-erdmann98}.
\end{proof}

\end{document}